\documentclass[11pt]{article}
\usepackage{amsfonts,amssymb,amsmath}

\usepackage[english]{babel}
\newtheorem{theorem}{Theorem}
\newtheorem{definition}{Definition}

\newtheorem{proof}{Proof}

\newtheorem{proposition}{Proposition}
\newtheorem{lemma}{Lemma}

\newcommand{\beq}{\begin{eqnarray}}
\newcommand{\eeq}{\end{eqnarray}}

\newcommand{\beqt}{\begin{eqnarray*}}
\newcommand{\eeqt}{\end{eqnarray*}}

\newcommand{\be}{\begin{equation}}
\newcommand{\ee}{\end{equation}}

\newcommand{\bl}{\begin{lemma}}
\newcommand{\el}{\end{lemma}}

\newcommand{\bt}{\begin{theorem}}
\newcommand{\et}{\end{theorem}}

\newcommand{\bd}{\begin{definition}}
\newcommand{\ed}{\end{definition}}

\newcommand{\bp}{\begin{proposition}}
\newcommand{\ep}{\end{proposition}}

\newcommand{\bpr}{\begin{proof}}
\newcommand{\epr}{\end{proof}}

\newcommand{\bi}{\begin{itemize}}
\newcommand{\ei}{\end{itemize}}

\newcommand{\ben}{\begin{enumerate}}
\newcommand{\een}{\end{enumerate}}

\newcommand{\Z}{\mathbb Z}
\newcommand{\R}{\mathbb R}
\newcommand{\N}{\mathbb N}

\newcommand{\E}{\mathbb E}

\newcommand{\pee}{\mathbb P}

\newcommand{\om}{\ensuremath{\omega}}
\newcommand{\Om}{\ensuremath{\Omega}}

\newcommand{\si}{\ensuremath{\sigma}}
\newcommand{\eps}{\ensuremath{\epsilon}}

\begin{document}

\title{{\bf  Transient random walks on 2d-oriented lattices}}

\author{Nadine Guillotin-Plantard\footnote{Universit\'e Claude Bernard - Lyon I, institut Camille Jordan,
b\^atiment Braconnier, 43 avenue du 11 novembre 1918,  69622
Villeurbanne Cedex, France. E-mail:
nadine.guillotin@univ-lyon1.fr}\footnote{Research partially
support by the RDSES program of the ESF.}, Arnaud Le
Ny\footnote{Universit\'e de Paris-Sud, laboratoire de
math\'ematiques, b\^atiment 425, 91405 Orsay cedex, France.
E-mail: arnaud.leny@math.u-psud.fr}}

\maketitle

\begin{center}
{\bf Abstract:}
\end{center}
We study the asymptotic behavior of the simple random walk on
oriented versions of $\Z^2$. The considered lattices are not
directed  on the vertical axis but unidirectional on the
horizontal one, with random orientations whose distributions are
generated by a dynamical system. We find a sufficient condition on
the smoothness of the generation for the transience of the simple
random walk on almost every such oriented lattices, and as an
illustration we provide a wide class of examples of inhomogeneous
or correlated distributions of the orientations. For ergodic
dynamical systems, we also prove a strong law of large numbers
and, in the particular case of i.i.d. orientations, we solve an
open problem and  prove a functional limit theorem in the space
$\mathcal{D}([0,\infty[,\R^2)$ of c\`adl\`ag functions, with an
unconventional normalization.
\medskip

%\footnotesize

\vspace{7cm}



\newpage
\section{Introduction}

The use  of random walks as a tool in
mathematical physics is now well established and they have been for
example widely used in classical statistical mechanics to study critical phenomena
 \cite{FFS}. Analogous methods in
  quantum statistical mechanics require the study of random walks on oriented
  lattices,
  due to the intrinsic non commutative character of the (quantum) world
\cite{CP2,LR}.
   Although random walks in random and non-random environments have been
   intensively studied for
    many years, only a few  results on random walks on oriented lattices are known. The
      recurrence or transience properties of simple random walks on oriented
versions of $\Z^2$ are studied in \cite{CP} when the horizontal
lines are
       unidirectional towards a random or deterministic direction. An interesting
       feature of this model is that, depending on the orientation, the walk could
       be either recurrent or transient. In a particular deterministic  case, alternatively  oriented horizontally rightwards
       and leftwards, the recurrence of the simple random walk is proved, whereas the transience
         naturally arises when the orientations are all identical in infinite regions.
         More surprisingly, it is also proved that the recurrent character of the simple
          random walk on $\Z^2$ is lost when the orientations are i.i.d. with zero mean.

In this paper, we study more general models and focus on spatially
inhomogeneous or dependent distributions  of the orientations. We
introduce lattices for which the distribution of the orientation
is generated by a dynamical system and prove that the transience
of the simple random walk still holds under smoothness conditions
on this generation. We detail examples and counterexamples for
various standard dynamical systems. For ergodic dynamical systems,
we also prove a strong law of large numbers and, in the case of
i.i.d. orientations, a functional limit theorem with an
unconventional normalization due to the random character of the
environment of the walk, solving an open question of \cite{CP}.

The model and our results are stated in Section 2, Section 3 is
devoted to the proofs, while illustrative examples of such
dynamical orientations are given in Section 4.

\section{Model and results}
\subsection{Dynamically oriented lattices}
Let $S=(E,{\cal A},\mu,T)$ be a dynamical system where $(E,{\cal
A},\mu)$ is a probability space and $T$ is an invertible
transformation of $E$ preserving the measure $\mu$. This system is
used to introduce inhomogeneity or dependencies in the
distribution of the random orientations, together with a function
$f$ from $E$ to $[0,1]$, which satisfies $\int_E f d\mu=1/2$ to
avoid trivialities. By {\em orientations}, we mean a random field
$\eps=(\eps_y)_{y \in \mathbb{Z}} \in \{-1,+1\}^{\mathbb{Z}}$,
i.e. a family $\eps$ of $\{-1,+1\}$-valued random variables
$\eps_y,\;y \in \Z$, and we distinguish two different approaches
to introduce its distribution.

\subsubsection{Quenched case:} It describes orientations spatially inhomogeneously distributed. For
$x \in E$ the {\em quenched law} $\mathbb{P}_{T}^{(x)}$ is the
product probability measure on $\{-1,+1\}^{\mathbb{Z}}$, equipped
with the product $\sigma-$algebra
$\mathcal{F}=\mathcal{P}(\{-1,+1\})^{\otimes \mathbb{Z}}$, whose
marginals can be  given by:
\[
\pee_{T}^{(x)}(\eps_y= + 1) =f(T^{y} x).
\]
To simplify, we have used the same notation for the quenched law
and its marginals, which should be written
$\mathbb{P}_{T,y}^{(x)}$ with $\mathbb{P}_{T}^{(x)}= \otimes_y
\mathbb{P}_{T,y}^{(x)}$. This quenched case is  an extension of
the i.i.d. case, with independent but not necessarily identically
distributed random variables. These random variables can be viewed
as the increments of a dynamic random walk \cite{gui0, gui1}.

\subsubsection{Annealed case:} We average on $x \in E$: the distribution of
$\eps$ is now  $\mathbb{P}_\mu$  defined for all $A \in
\mathcal{F}$ by

\[
\mathbb{P}_\mu[\eps \in A] = \int_E \pee_T^{(x)} [\eps \in A] d
\mu(x).
\]
The  marginals are thus given for all $y \in \mathbb{Z}$ by
\[
\mathbb{P}_\mu[\eps_y=+1] = \int_E f(T^yx)  d \mu(x)=\int_E f(x)
d\mu(x)=\frac{1}{2}
\]
and the hypothesis $\int_E f d\mu =\frac{1}{2}$ has been taken to
get $\mathbb{E}_\mu[\eps_y]=0$. The $T$-invariance of $\mu$
implies the translation-invariance of $\mathbb{P}_\mu$ but this
latter is  not a product measure in general: The correlations of
the dynamical system for the function $f$, defined  for all $y \in
\Z$ by
\begin{eqnarray}\label{cmu} C_\mu^f(y)&:=&\int_E f(x) \cdot f
\circ T^y(x) d\mu(x) - \int_E f(x)
d\mu(x) \cdot\int_E f \circ T^y(x) d \mu(x) \\
&=& \int_E f(x) f(T^y(x)) d\mu(x) - \frac{1}{4} \nonumber
\end{eqnarray}
 are indeed related with the covariance of the $\eps$'s after a short
 computation:
\begin{equation}\label{correlation} \forall y \in \Z, \rm{Cov}_\mu (\eps_0, \eps_y) =
4 \; C_\mu^f(y). \end{equation} One can thus construct dependent
variables whose dependence is directly related to the correlations
of the dynamical system. This annealed case leads in Section 4 to
another extension of the i.i.d. case where, independence is
dropped but translation-invariance is kept.

\subsubsection{Lattices}

We use these orientations to build  {\em dynamically oriented
lattices}. They are oriented versions of $\Z^2$: the vertical
lines are not oriented and the horizontal ones are unidirectional,
the orientation at a level $y \in \Z$ being given by the random
variable $\eps_y$ (say right if the value is $+1$ and  left if
 it is $-1$). More formally we give the
\bd[Dynamically oriented lattices] Let $\eps=(\eps_y)_{y \in \Z}$
be a sequence of orientations defined as previously. The {\em
dynamically oriented lattice}
$\mathbb{L}^\eps=(\mathbb{V},\mathbb{A}^\eps)$ is the (random)
directed graph with (deterministic) vertex set $\mathbb{V}=\Z^2$
and (random) edge set $\mathbb{A}^\eps$ defined by the condition
that for $u=(u_1,u_2), v=(v_1,v_2) \in \Z^2$, $(u,v) \in
\mathbb{A}^\eps$ if and only if  $v_1=u_1$ and $v_2=u_2 \pm 1$, or
$v_2=u_2$ and $v_1=u_1+ \eps_{u_2}$.\ed

\subsection{Simple random walk on $\mathbb{L}^\eps$}

We consider the  simple random walk  $M=(M_n)_{n \in \mathbb{N}}$
on $\mathbb{L}^\eps$.
 For every realization $\eps$, it is a $\mathbb{Z}^2$-valued Markov chain defined on a probability space $(\Omega,
 \mathcal{B},\mathbb{P})$, whose ($\eps$-dependent) transition probabilities are
 defined for all $(u,v) \in \mathbb{V}\times \mathbb{V}$ by
\[\pee[M_{n+1}=v  | M_n=u]=\;\left\{
\begin{array}{lll} \frac{1}{3}  \; &\rm{if} \;  (u,v)   \in \mathbb{A}^\eps&\\
\\
0 \; \; &\rm{otherwise.}&
\end{array}
\right.
\]
Its transience is proved in \cite{CP} for almost every orientation
$\eps$ when the $\eps_y$'s are i.i.d. and  we generalize it in
this dynamical context when the orientations are either {\em
annealed} or {\em quenched}.

\bt \label{thm1} Assume that
\begin{equation} \label{C}
\int_{E}
\frac{1}{\sqrt{f(1-f)}}\ d\mu <\infty \end{equation}
 then:
\begin{enumerate} \item In the annealed case, for
$\mathbb{P}_{\mu}$-a.e. orientation $\eps$, the simple random walk
on dynamically oriented lattice $\mathbb{L}^{\eps}$ is transient.
\item In the quenched case, for $\mu$-a.e.  $x\in E$, for
$\mathbb{P}_{T}^{(x)}$-a.e. realization of the orientation $\eps$,
the simple random walk on the dynamically oriented lattice
$\mathbb{L}^\eps$ is transient.
\end{enumerate}
\et {\bf Remarks}
\begin{enumerate}
\item Non-invertible transformations $T$ of the space $E$ can also
be considered and in this case it is straightforward to extend the
conclusions of  Theorem \ref{thm1} if the distribution of the
orientations $(\epsilon_{y})_{y\in \Z}$ have  marginals defined by
$\pee_{T}^{(x)} (\epsilon_{y} = +1) = f(T^{|y|} x)$. The measure
$\pee_{\mu}$ is not stationary anymore in the annealed case (see
the example of Manneville-Pomeau maps of the interval in Section
4). \item In Section 4 we exhibit dynamical systems for which
Theorem \ref{thm1} applies (Bernoulli or Markov Shifts,
Manneville-Pomeau maps), but also counter-examples from a family
of dynamical systems (irrational and rational rotations on the
torus). The latter case provides instructive examples: when the
function $f$ satisfies (\ref{C}), the ergodicity (or not) of the
dynamical system is not required, whereas when $f$ does not fulfil
(\ref{C}), the properties of the underlying dynamical system can
play a role, e.g. in the non-ergodic case when, according to the
rational angle we choose, the simple random walk on the
corresponding oriented lattice can be transient or recurrent.
\end{enumerate}
\subsection{Limit theorems in the ergodic case}
 Let us assume that the
dynamical system $S=(E,{\cal A},\mu,T)$ defined in Section 2.1 is
ergodic. \bt[Strong law of large numbers] The random walk on the
lattice $\mathbb{L}^{\eps}$ has $\pee \otimes \pee_\mu$-almost
surely zero speed, i.e.
\begin{equation}
\lim_{n\rightarrow +\infty}\frac{M_{n}}{n}=(0,0)\ \ \ \ \ \pee
\otimes \pee_\mu-{\rm almost} \; {\rm surely.}
\end{equation}
\et

\subsubsection{Functional limit theorem for i.i.d. orientations}

We also answer in this paper to an open question of \cite{CP} and
obtain a functional limit theorem with a suitable  normalization.
We establish that the study of the simple random walk on
$\mathbb{L}^\eps$ is closely related to a {\it simple random walk
in a random scenery} defined for every $n\ge 1$ by
$$Z_n=\sum_{k=0}^{n} \eps_{Y_k}$$
where $(Y_k)_{k\ge 0}$ is the simple random walk on $\Z$ starting
from 0. Consider a standard Brownian motion $(B_{t})_{t\ge 0}$,
denote by $(L_{t}(x))_{t \ge 0}$ its corresponding local time at
$x\in\R$ and introduce a pair of independent Brownian motions
$(Z_{+}(x), Z_{-}(x)), x\geq 0$ defined on the same probability
space as $(B_{t})_{t\ge 0}$ and independent of him. The following
process is well-defined for all $t \geq 0$:
\begin{equation}\label{th}
\Delta_{t}=\int_{0}^{\infty}L_{t}(x)dZ_{+}(x)+\int_{0}^{\infty}L_{t}(-x)dZ_{-}(x).
\end{equation}
It has been proved by Kesten {\em et al.} \cite{KS} that this
process has a self-similar continuous version  of index
$\frac{3}{4}$, with stationary increments. We denote
$\stackrel{\mathcal{D}}{\Longrightarrow}$ for a  convergence in
the space of c\`adl\`ag functions $\mathcal{D}([0,\infty),\R)$
endowed with the Skorohod topology.
 \bt\label{thm11} \mbox{\bf [Kesten and Spitzer (1979)]}
 \be \Big(\frac{1}{n^{3/4}} Z_{[nt]} \Big)_{t \geq 0}
\; \stackrel{\mathcal{D}}{\Longrightarrow} (\Delta_t)_{t \geq 0}.
\ee  \et

We introduce a real constant $m=\frac{1}{2}$, defined later as the
mean of some geometric random variables related to the behavior of
the walk in the horizontal direction\footnote{Our results are in
fact valid for similar model for which $m \neq \frac{1}{2}$
corresponding to non symmetric nearest neighbors random walks. Of
course the transience is not at all surprising in this case, but
getting the limit theorems can be of interest.}. Using Theorem
\ref{thm11}, we shall prove

\bt[Functional limit theorem]\label{thm2}
 \be \label{flt}
\Big(\frac{1}{n^{3/4}} M_{[nt]} \Big)_{t \geq 0} \;
\stackrel{\mathcal{D}}{\Longrightarrow} \frac{m}{(1+m)^{3/4}}(
\Delta_t,0)_{t \geq 0}. \ee \et {\bf Remark} : It is not
surprising that the vertical component is negligible towards
$n^{3/4}$ because its fluctuations are of  order $\sqrt{n}$. We
suspect that we have in fact
\[
\Big(\frac{1}{n^{3/4}} M^{(1)}_{[nt]}, \frac{1}{n^{1/2}
}M^{(2)}_{[nt]}\Big)_{t \geq 0} \;
\stackrel{\mathcal{D}}{\Longrightarrow} \Big(\frac{m}{(1+m)^{3/4}}
\Delta_t,B_t \Big)_{t \geq 0}
\]
but this is not straightforward because  the horizontal
($M^{(1)}$) and vertical ($M^{(2)}$)  components are not
independent. We  believe that  $(B_t)_{t \geq 0}$ and
$(\Delta_t)_{t \geq 0}$ are independent but this also has to be
proved.

%We conjecture the following local limit theorem for the random
%walk $M_n$.
% \bcon[Local limit theorem]\label{conj} There exists a constant $C>0$ such
%that as $n\rightarrow +\infty$,
 %\be \label{tll}
%\tilde{\pee}_\mu [M_{n}=(0,0) ]\sim C n^{-5/4}. \ee \econ

\section{Proofs}
\subsection{Vertical and horizontal embeddings of the simple random walk}
 The simple random walk $M$  defined  on $(\Om,\mathcal{B},\pee)$
 can be decomposed into vertical and horizontal embeddings by projection to
 the corresponding  axis. The vertical one is a simple random walk
 $Y=(Y_n)_{n \in \N}$ on $\mathbb{Z}$ and we define for all $n \in \N$
 its {\em local time} at the level $y \in \Z$ by
\[
\eta_n(y)=\sum_{k=0}^n \mathbf{1}_{Y_k=y}.
\]
The horizontal embedding is a random walk with $\mathbb{N}$-valued
geometric jumps: a doubly infinite family $(\xi_i^{(y)})_{i \in
\mathbb{N}^*, y \in \Z}$ of independent geometric random variables
of mean $m=\frac{1}{2}$ is given and one defines the embedded
horizontal random walk $X=(X_n)_{n\in\N}$ by $X_0=0$ and for $n
\geq 1$,
\[
X_n=\sum_{y \in \Z} \eps_y \sum_{i=1}^{\eta_{n-1}(y)} \xi_i^{(y)}
\]
with the convention that the last sum is zero when
$\eta_{n-1}(y)=0$. Of course, the walk $M_n$ does not coincide
with $(X_n,Y_n)$ but these objects are closely related: Define for
all $n \in \N$
\[
T_n=n + \sum_{y \in \Z} \sum_{i=1}^{\eta_{n-1}(y)} \xi_i^{(y)}
\]
to be the instant just after the random walk $M$ has performed its
n$^{\rm{th}}$ vertical move. A direct and useful consequence of
this decomposition is the following result \cite{CP}.
 \bl \label{lem1} \ben
\item $M_{T_n}=(X_n,Y_n),\; \forall n \in \N$. \item For a given
orientation $\eps$, the transience of $(M_{T_n})_{n \in \N}$
implies the transience of $(M_n)_{n \in \N}$. \een \el
\subsection{Proof of the transience of the simple random walk}

The vertical walk $Y$, independent of $\eps$, is known to be
recurrent with fluctuations of order $\sqrt{n}$.  For any $i \in
\N$, $\delta_i$ is  a strictly positive  real number and we write
$d_{n,i}=n^{\frac{1}{2}+\delta_i}$ to introduce a partition of
$\Omega$ between typical or untypical paths of $Y$:
\[
A_n=\big\{ \om \in \Om; \max_{0 \leq k \leq 2n} \; |Y_k|  <
d_{n,1} \big\} \; \cap \; \big\{ \om \in \Om; \max_{y \in \Z} \;
\eta_{2n-1}(y) < d_{n,2}\big\}
\]
and
\[
B_n=\big\{\om \in A_n; \Big| \sum_{y \in \Z} \eps_y \eta_{2n-1}(y)
\Big|  > d_{n,3}\big\}.
\]
We first consider the joint measures $\tilde{\mathbb{P}}_\mu =
\pee \otimes \pee_\mu$ (annealed case) or $\tilde{\pee}_{T}^{(x)}=
\pee \otimes \pee_{T}^{(x)}$ (quenched case) and prove that \be
\label{eqn1} \sum_{n \in \N} \tilde{\mathbb{P}}_\mu
[X_{2n}=0;Y_{2n}=0] \; < \; \infty.
 \ee

%\begin{equation}\label{eqn0}
% \tilde{\mathbb{P}}_\mu [M_{T_n}=(0,0) \; {\rm i.o}.] =0
%\end{equation}
%or, using Lemma \ref{lem1},

By definition
\[
 \sum_{n \in \N} \tilde{\mathbb{P}}_\mu [X_{2n}=0;Y_{2n}=0] =
 \int_E\sum_n \pee[ \pee_{T}^{(x)}[X_{2n}=0;Y_{2n}=0] ]d \mu(x)
\]
and we first decompose $\tilde{\pee}_{T}^{(x)}[X_{2n}=0;Y_{2n}=0]$
into \[ \tilde{\pee}_{T}^{(x)}[X_{2n}=0;Y_{2n}=0;A_n^c] +
\tilde{\pee}_{T}^{(x)}[X_{2n}=0;Y_{2n}=0;B_n]
 + \tilde{\pee}_{T}^{(x)}[X_{2n}=0;Y_{2n}=0; A_n \setminus B_n].
\]
 Some results of the i.i.d. case of \cite{CP} still hold uniformly in $x$ and in particular
 we can prove using standard techniques the following
\bl \label{lem3} \begin{enumerate} \item For every $x\in E$,
$\sum_{n \in \N} \tilde{\pee}_{T}^{(x)}[X_{2n}=0;Y_{2n}=0;A_n^c]
\; < \; \infty$. \item For every $x\in E$, $\sum_{n \in \N}
\tilde{\pee}_{T}^{(x)}[X_{2n}=0;Y_{2n}=0;B_n] \; < \; \infty$.
\end{enumerate}
\el Define the $\si$-algebras  $\mathcal{F}=\sigma(Y)$ and
$\mathcal{G}=\sigma(\eps)$  generated by the families of r.v.'s
$Y$ and $\eps$. Then one has
$$
p_n^{(x)}: =\tilde{\pee}_{T}^{(x)}[X_{2n}=0;Y_{2n}=0;A_n \setminus
B_n] = \E\Big[ \mathbf{1}_{Y_{2n=0}} \E \big[\mathbf{1}_{A_n
\setminus B_n} \tilde{\pee}_{T}^{(x)} \big[X_{2n}=0 \big|
\mathcal{F} \vee \mathcal{G} \big] \big| \mathcal{F} \big] \Big].
$$.
To prove the theorem, it remains to show that \be \label{pn}\int_E
\left(\sum_{n \in \N} p_n^{(x)}\right)\ d\mu(x)<\infty. \ee Recall
that for the simple random walk $Y$, there exists $C>0$ s.t. \be
\label{srw} \pee[Y_{2n}=0] \sim C \cdot n^{-\frac{1}{2}}, \; n
\rightarrow + \infty \ee and we can prove as in \cite{CP} the

\bl \label{lem5}

On the set $A_n \setminus B_n$, we have uniformly in $x\in E$, \be
\label{sqrtln} \tilde{\pee}_{T}^{(x)} \big[X_{2n}=0 \big|
\mathcal{F} \vee \mathcal{G}\big]= \mathcal{O} \Big(
\sqrt{\frac{\ln{n}}{n}}\Big). \ee

\el Hence, the transience of the simple random walk is a direct
consequence of the following

\bp \label{prop1}
It is possible to choose
$\delta_1,\delta_2,\delta_3>0$ such that there exists $\delta>0$
and
 \be \label{eqn3}
\int_E \tilde{\pee}_{T}^{(x)}\big[ A_n \setminus B_n \big|
\mathcal{F} \big]\ d\mu(x) = \mathcal{O} \big(n^{-\delta}). \ee
\ep

{\bf Proof :}  We have to estimate, on the event $A_{n}$, the
conditional probability
$$\tilde{\pee}_{T}^{(x)}[|\sum_{y\in\mathbb{Z}}\zeta_{y}|\le
d_{n,3} \big| \mathcal{F} \big]$$
 where
$\zeta_{y}=\eps_y \eta_{2n-1}(y), y\in\mathbb{Z}$. Let $G$ be a
centered Gaussian random variable with variance $d_{n,3}^2$,
(conditionally on $\mathcal{F}$) independent of the random
variables $\zeta_y$'s. Clearly,
$$\tilde{\pee}_{T}^{(x)}\big[\sum_{y}\zeta_y\in [0,d_{n,3}]\big| \mathcal{F}\big] =
\frac{\tilde{\pee}_{T}^{(x)}\big[\sum_{y}\zeta_y\in [0,d_{n,3}] ;
0\le G\le d_{n,3}\big| \mathcal{F}\big]}
{\tilde{\pee}_{T}^{(x)}[0\le G\le d_{n,3}\big| \mathcal{F}]}$$
where $\tilde{\pee}_{T}^{(x)}[0\le G\le d_{n,3}\big|
\mathcal{F}]=c>0$ is independent of $n$. Since $G$ is independent
of the random variables $\zeta_{y}$'s and using the symmetry of
the Gaussian distribution, we have
$$
\tilde{\pee}_{T}^{(x)}\big[\sum_{y}\zeta_y\in [0,d_{n,3}] ; 0\le
G\le d_{n,3}\big|
\mathcal{F}\big]=\tilde{\pee}_{T}^{(x)}\big[\sum_{y}\zeta_y\in
[0,d_{n,3}] ; -d_{n,3}\le G\le 0\big| \mathcal{F}\big].$$
Consequently, we obtain
$$\tilde{\pee}_{T}^{(x)}\big[\sum_{y}\zeta_y\in [0,d_{n,3}]\big| \mathcal{F}]
\le \frac{1}{c}\tilde{\pee}_{T}^{(x)}[|\sum_{y}\zeta_y +G|\le
d_{n,3} \big| \mathcal{F}\big] \ \rm{and}$$
$$\tilde{\pee}_{T}^{(x)}\big[\sum_{y}\zeta_y\in [-d_{n,3},0]\big|
\mathcal{F}\big]\le
\frac{1}{c}\tilde{\pee}_{T}^{(x)}\big[|\sum_{y}\zeta_y +G|\le
d_{n,3}\big| \mathcal{F}\big]$$ and then, we have the following
inequality
$$\tilde{\pee}_{T}^{(x)}\big[|\sum_{y}\zeta_{y}|\le d_{n,3}\big|
\mathcal{F}\big]\le \frac{2}{c}\
\tilde{\pee}_{T}^{(x)}\big[|\sum_{y}\zeta_y +G|\le d_{n,3}\big|
\mathcal{F}\big].$$ From Plancherel's formula, we deduce that
there exists a constant $C>0$ such that \be \label{eqn4}
\tilde{\pee}_{T}^{(x)}\big[ |\sum_{y}\zeta_y +G|\le d_{n,3} \big|
\mathcal{F} \big] \leq C \cdot d_{n,3} \cdot I_n(x) \ee where
\[
I_n(x)=\int_{-\pi}^{\pi} \E\big[e^{it \sum_{y \in \Z} \eps_y
\eta_{2n-1}(y)} \big| \mathcal{F}  \big]e^{-t^2 d_{n,3}^2/2} dt.
\]
To use that for $td_{n,3}$ small enough, $e^{-t^2 d_{n,3}^2/2}$
dominates the term under the expectation, we split the integral in
two parts. For $b_n=\frac{n^{\delta_2}}{d_{n,3}}$, we write
$I_n(x)=I_n^{1}(x) + I_n^2(x)$ with \beqt I_n^1(x)=\int_{|t| \leq
b_n} \E\big[e^{it \sum_{y \in
\Z} \eps_y \eta_{2n-1}(y)} \big| \mathcal{F} \big]  e^{-t^2 d_{n,3}^2/2} dt\\
I_n^2(x)=\int_{|t| > b_n} \E \big[e^{it \sum_{y \in \Z} \eps_y
\eta_{2n-1}(y)} \big| \mathcal{F} \big]e^{-t^2 d_{n,3}^2/2} dt.
\eeqt To control the integral $I_n^2(x)$, we write \beqt
|I_n^2(x)| &\leq& C \int_{|t| > b_n} e^{-t^2 d_{n,3}^2/2} dt=
\frac{C}{d_{n,3}} \int_{|s| > n^{\delta_2}} e^{-s^2/2} ds \;  \leq
\;  \frac{2C}{d_{n,3}} \/  n^{-\delta_2} \/  e^{-n^{2 \delta_2}
/2} \eeqt to get uniformly in $x\in E$
\[
|I_n^2(x)|=\mathcal{O} \big(e^{-n^{2 \delta_2} / 2}).
\]

\bl \label{lem6} For $\delta_3 > 2 \delta_2 $,
\[
\int_E |I_n^1(x)|\ d\mu(x)=\mathcal{O}
\big(n^{-\frac{3}{4}+\frac{\delta_1}{2}}\big).
\]
\el

{\bf Proof :} From the definition of the orientations
$(\eps_{y})_{y}$, an explicit formula for the characteristic
function $\phi_{\eps_y}^{(x)}$ of the random variable $\eps_y$ can
be given and we deduce that
 \[
|\phi_{\eps_y}^{(x)}(u)|^2 =\cos^2(u)+(2f(T^{y} x)-1)^2\sin^2(u) =
1-4f(T^{y} x)(1-f(T^{y} x))\sin^2(u) \] and by independence of the
$\eps$'s we get
$$|I_n^1(x)|\le \int_{|t| \leq b_n} |
\prod_{y}\phi_{\eps_y}^{(x)}(\eta_{2n-1}(y)t)|\ dt.$$ Denote
$p_{n,y}=\frac{\eta_{2n-1}(y)}{2n}$, $C_n=\{y:\eta_{2n-1}(y)\ne
0\}$ and use H\"older's inequality to get \[ |I_n^1(x)|\leq
\prod_{y} \Big[ \Big( \int_{|t| \leq b_n}
|\phi_{\eps_y}^{(x)}(\eta_{2n-1}(y)t)|^{1/p_{n,y}} dt
\Big)^{p_{n,y}}\Big]. \] Now, using the fact that we work on
$A_n$, we choose $\delta_3>2 \delta_2$ s.t.
 $\lim_n b_n\eta_{2n-1}(y)= 0$ uniformly in $y$.
Using  $\sin(x)\ge \frac{2}{\pi} x$ for $x\in [0,\frac{\pi}{2}]$
and $\exp(-x)\ge 1-x$, one has \begin{eqnarray*} |I_n^1(x)|& \leq
& \prod_{y\in C_n} \left(\frac{1}{\eta_{2n-1}(y)} \int_{|v|\leq
b_n \eta_{2n-1}(y)}\exp\left(-\frac{16}{p_{n,y}\pi^2}f(T^{y}
x)(1-f(T^{y}x))v^2\right)\ dv \right)^{p_{n,y}}\\
& \leq & \prod_{y\in C_n} \left(\frac{c
\textbf{1}_{f(T^{y}x)(1-f(T^{y}x))>0}}{\sqrt{2n\eta_{2n-1}(y)f(T^{y}x)(1-f(T^{y}x))}
}\right)^{p_{n,y}} \ \ (\mbox{with }\ \ c=\pi^{3/2}/{4})\\
& = & c  \exp\big[-\frac{1}{2} \sum_{y\in C_n}
p_{n,y}\log(2n\eta_{2n-1}(y))\big]\cdot \prod_{y\in C_n}
\left(\frac{\textbf{1}_{f(T^{y}x)(1-f(T^{y}x))>0}}{\sqrt{f(T^{y}x)(1-f(T^{y}x))}
}\right)^{p_{n,y}}. \end{eqnarray*} The vector
$\textbf{p}=(p_{n,y})_{y \in C_n}$ defines a probability measure
on  $C_n$ and we have \[ -\frac{1}{2}\sum_{y\in C_n}
p_{n,y}\log(2n\eta_{2n-1}(y))=  -\log 2n -\frac{1}{2}\sum_{y\in
C_n} p_{n,y}\log p_{n,y}  =  -\log 2n+ \frac{1}{2}H(\textbf{p})
\]
 where $H(\cdot)$ is the entropy of the probability vector
$\textbf{p}$, always bounded by $\log(\textrm{card}(C_n))$. We
thus have on the set $A_n$,
\[
|I_n^1(x)| \leq  c\exp \left[-\log 2n + \frac{1}{2}\log (2 d_{n,1})\right]
\prod_{y\in C_n} \left(\frac{\textbf{1}_{\{f(T^{y}x)(1-f(T^{y}x))>0\}}}{\sqrt{f(T^{y}x)(1-f(T^{y}x))}}\right)^{p_{n,y}}.\]
 By applying H\"{o}lder's inequality and the fact that $T$ preserves the measure $\mu$, we get
 \beqt
\int_{E}|I_n^1(x)|\ d\mu(x)& \leq & C\cdot n^{-\frac{3}{4}+\frac{\delta_{1}}{2}}
\int_{E}\prod_{y\in C_n}\left(\frac{\textbf{1}_{\{f(T^{y}x)(1-f(T^{y}x))>0\}}}{\sqrt{f(T^{y}x)(1-f(T^{y}x))}
}\right)^{p_{n,y}}\ d\mu(x)\\
&\leq & C \cdot n^{-\frac{3}{4}+\frac{\delta_{1}}{2}}
\prod_{y\in C_n}\left[\int_{E}\left(\frac{\textbf{1}_{\{f(T^{y}x)(1-f(T^{y}x))>0\}}}{\sqrt{f(T^{y}x)(1-f(T^{y}x))}
}\right)\ d\mu(x)\right]^{p_{n,y}}\\
&=& C\cdot  n^{-\frac{3}{4}+\frac{\delta_{1}}{2}}
\int_{E}\frac{1}{\sqrt{f(x)(1-f(x))}}\ d\mu(x).\;  \;  \;  \;  \;
\;  \;  \;  \; \;  \;  \;  \; \;  \;  \;  \; \;  \;  \;  \; \;  \;
\;  \; \;  \;  \;  \; \;  \;  \;  \; \;  \;
 \diamond
 \eeqt

Now, using (\ref{eqn4}), write with the usual notation
$d_{n,3}= n^{\frac{1}{2}+\delta_3}$:
\[
\int_{E} \tilde{\pee}_{T}^{(x)}[A_n \setminus B_n | \mathcal{F}]\
d\mu(x) \leq C\cdot d_{n,3} \int_{E}\left(|I_n^1(x)|+
|I_n^2(x)|\right)\ d\mu(x)
\]
 and consider $\delta_3 > 2 \delta_2$. By the previous lemmata,   we
have
\[
d_{n,3}\cdot \int_{E}|I_n^1(x)|\ d\mu(x)= \mathcal{O}
\big(n^{-\frac{1}{4}+\delta_3+\frac{\delta_1}{2}}\big), \; d_{n,3}
\cdot \int_{E}|I_n^2(x)|\ d\mu(x)= \mathcal{O} \big(e^{-n^{2
\delta_2} / 2})
\]
and the proposition follows by choosing $\delta_1, \delta_2,
 \delta_3$ small enough. $\;  \;  \; \;  \;  \;  \; \;  \;\;  \;  \;
  \;  \;  \;  \; \;  \;\;  \;  \; \;  \;  \;  \; \;  \;\;  \;  \; \;
   \;  \;  \;  \;  \; \;  \; \diamond$

Combining Equations  (\ref{srw}), (\ref{sqrtln}) and (\ref{eqn3}),
we obtain (\ref{pn}) and then (\ref{eqn1}). By Borel-Cantelli's
Lemma, we get :
\[
\tilde{\pee}_\mu \big[ M_{T_n}=(0,0)\;\rm{i.o.} \big]=\pee_\mu
\big[ \pee \big[ M_{T_n}=(0,0)\; \rm{i.o.} \big] \big] =0
\]
and thus for $\pee_\mu$-almost every orientation $\eps$, $\pee
\big[ M_{T_n}=(0,0) \; \rm{i.o.} \big] =0$. This proves that
$(M_{T_n})_{n \in \N}$ is transient for $\pee_{\mu}$-almost every
orientation $\eps$, and by Lemma \ref{lem1}, the $\pee_\mu$-almost
sure transience of the simple random walk on the annealed oriented
lattice. Transience in the quenched case is a direct consequence
of the transience in the annealed case.
\subsection{Proof of the strong law of large numbers}

\bl \label{slln}\mbox{\bf (SLLN for the embedded random walk)}
\begin{equation}
\lim_{n\rightarrow +\infty}\frac{M_{T_n}}{n}=(0,0)\ \
\mbox{$\tilde{\pee}_\mu$-almost surely.}
\end{equation}
\el

{\bf Proof :} Since $(Y_{n})_{n\geq 0}$ is a simple random walk,
$\frac{Y_{n}}{n}$ goes to $0$ as $n \to \infty$
$\tilde{\pee}_\mu$-a.s. as $n \to \infty$ and it is enough to
prove that $(\frac{X_{n}}{n})$ converges almost surely to 0.
Introduce
$$Z_{n}=\sum_{k=0}^{n-1} \eps_{Y_{k}}=\sum_{y\in\Z}\eps_{y} \eta_{n-1}(y).$$
Under the probability measure $\tilde{\pee}_{\mu}^{~}$, the
stationary sequence $(\eps_{Y_{k}})_{k\geq 0}$ is ergodic
\cite{Ka}, so from Birkhoff's theorem, as $n$ tends to infinity,
\[
\frac{Z_{n}}{n}\rightarrow \E[\eps_{0}]=0\ \ \mbox{almost surely.}
\]
Clearly, $X_{n}-mZ_{n}=\sum_{y\in\Z}
\eps_{y}\sum_{i=1}^{\eta_{n-1}(y)} (\xi_{i}^{(y)}-m)$ and for an
even integer $r$

\small
$$\E[(X_{n}-mZ_{n})^{r}]=\sum_{y_1\in\Z,\ldots y_{r}\in\Z} \E\left[\eps_{y_1}\ldots
\eps_{y_{r}}\sum_{i_{1}=1}^{\eta_{n-1}(y_{1})}\ldots
\sum_{i_{r}=1}^{\eta_{n-1}(y_{r})} \E[(\xi_{i_1}^{(y_1)}-m)\ldots
(\xi_{i_{r}}^{(y_{r})}-m)|\mathcal{F} \vee \mathcal{G}]\right].$$
\normalsize
 The $\xi_i^{(y)}$'s are independent of the vertical
walk and the orientations; moreover, the random variables
$\xi_{i}^{(y)}-m, i\geq 1, y\in \Z $ are i.i.d. and  centered, so
the summands are non zero if and only if $i_{1}=\ldots=i_{r}$ and
$y_{1}=\ldots =y_{r}$. Then,
$$\E[(X_{n}-mZ_{n})^{r}]=n\E[(\xi_{1}^{(0)} -m)^{r}]:=nm_{r}\; \; \; \; \; \mbox{   (say)}.$$
Let $\delta>0$. By Tchebychev's inequality,
\begin{eqnarray*}
\pee\Big[\left|\frac{X_{n}-mZ_{n}}{n}\right|\geq
\epsilon\Big]&\leq &\frac{1}{\delta^{r} n^{r}}
\E[(X_{n}-mZ_{n})^{r}] \; \leq \; \frac{m_{r}}{\delta^{r}
n^{r-1}}.
\end{eqnarray*}
We choose $r=4$ and thus from Borel-Cantelli Lemma, we deduce that
$\frac{X_{n}-mZ_{n}}{n}$ converges almost surely to 0 as $n$ goes
to infinity. $\diamond$

Using similar techniques, one also proves the

 \bl \label{tn} The sequence  $(\frac{T_{n}}{n})_{n\ge 1}$
converges $\tilde{\pee}_\mu$-a.s. to $(1+m)$ as $n\rightarrow
~+\infty.$ \el

%{\bf Proof :}  Remark first that
%$T_{n}=n+\sum_{y\in\Z}\sum_{i=1}^{\eta_{n-1}(y)}(\xi_{i}^{(y)}-m)
%+m\sum_{y\in\Z}\eta_{n-1}(y).$ Now, \beqt
%\E\Big[\Big(\sum_{y\in\Z}\sum_{i=1}^{\eta_{n-1}(y)}(\xi_{i}^{(y)}-m)\Big)^3\Big]
%&=&\sum_{y_{1}\in\Z,\ldots,y_{3}\in\Z}
%\E\Big[\sum_{i_{1}=1}^{\eta_{n-1}(y_{1})}\ldots
%\sum_{i_{3}=1}^{\eta_{n-1}(y_{3})}(\xi_{i_{1}}^{(y_{1})}-m)
%\ldots (\xi_{i_{3}}^{(y_3)}-m)\Big]\\
%&=& \sum_{y_{1}\in\Z,\ldots,y_{3}\in\Z}
%\E\Big[\sum_{i_{1}=1}^{\eta_{n-1}(y_{1})} \ldots
%\sum_{i_{3}=1}^{\eta_{n-1}(y_{3})}\E[(\xi_{i_{1}}^{(y_{1})}-m)\ldots
%(\xi_{i_{3}}^{(y_{3})}-m)|\mathcal{F} \vee
%\mathcal{G}]\Big]\\
%&=& m_{3} \sum_{y_{1}\in\Z,\ldots,y_{3}\in\Z}
%\E\Big[\sum_{i_{1}=1}^{\eta_{n-1}(y_{1})}
%\ldots \sum_{i_{3}=1}^{\eta_{n-1}(y_{3})} {\bf 1}_{i_{1}=\ldots =i_{3}}{\bf 1}_{y_{1}=\ldots =y_{3}}\Big]\\
%&=& m_3 \sum_{y\in\Z} \E[\eta_{n-1}(y)]= m_3 n\ \ (\mbox{\rm
%where\ } \ m_{3}=\E[(\xi_{1}^{(0)}-m)^3]). \eeqt

%From Tchebychev's inequality and Borel-Cantelli Lemma again, we
%deduce that, as $n\rightarrow ~+\infty$,
%$$\frac{1}{n}\sum_{y\in\Z}\sum_{i=1}^{\eta_{n-1}(y)}(\xi_{i}^{(y)}-m)\rightarrow
%0 \ \mbox{\rm a.s.}.$$ Using the fact that
%$\sum_{y\in\Z}\eta_{n-1}(y) =n$, we deduce the lemma. $\diamond$

Let us prove now the almost sure convergence of the sequence
$(\frac{M_{n}}{n})_{n\geq 1}$ to $(0,0)$. Since the sequence
$(T_{n})_{n\geq 1}$ is strictly increasing, there exists a
non-decreasing sequence of integers sequence $(U_{n \geq 1})_n$
such that $T_{U_{n}}\leq n< T_{U_{n}+1}$. Denote
$M_n=(M_n^{(1)},M_n^{(2)})$, then we have $ M_{n}^{(1)} \in
[\min(M_{T_{U_{n}}}^{(1)},M_{T_{U_{n}+1}}^{(1)}),\max(M_{T_{U_{n}}}^{(1)},M_{T_{U_{n}+1}}^{(1)})]$
and $M_{n}^{(2)}=M_{T_{U_{n}}}^{(2)}$, by definition of the
embedding. The (sub-)sequence $(U_{n})_{n\geq 1}$ is nondecreasing
and $\lim_{n\rightarrow +\infty}U_n=+\infty$, and then by
combining Lemmata \ref{slln} and \ref{tn}, we get that as
$n\rightarrow +\infty$, \be\label{tun}
\frac{M_{T_{U_{n}}}}{T_{U_{n}}}\rightarrow (0,0) \mbox{
$\tilde{\pee}_\mu$ a.s.} \ee Now,
\[
\displaystyle\left|\frac{M_{n}^{(1)}}{n}\right|\leq
\max\left(\left|\frac{M_{T_{U_{n}}}^{(1)}}{n}\right|,
\left|\frac{M_{T_{U_{n}+1}}^{(1)}}{n}\right|\right)\leq
\max\left(\left|\frac{M_{T_{U_{n}}}^{(1)}}{T_{U_{n}}}\right|,
\left|\frac{M_{T_{U_{n}+1}}^{(1)}}{T_{U_{n}}}\right|\right)
\]
and
$$\left|\frac{M_{n}^{(2)}}{n}\right|=
\left|\frac{M_{T_{U_{n}}}^{(2)}}{T_{U_{n}}}\right| \cdot
\frac{T_{U_{n}}}{n} \leq
\left|\frac{M_{T_{U_{n}}}^{(2)}}{T_{U_{n}}}\right|.$$ From
(\ref{tun}), we deduce the almost sure convergence of the
coordinates to 0 and then this of the sequence
$(\frac{M_{n}}{n})_{n\geq 1}$ to (0,0) as $n\rightarrow\infty$.

\subsection{Proof of the functional limit theorem}

\bp\label{pr1} The sequence of random processes
$n^{-3/4}(X_{[nt]})_{t\ge 0}$ weakly converges in the space ${\cal
D}([0,\infty[,\R)$ to the process $(m\Delta_{t})_{t\ge 0}$. \ep

{\bf Proof :} Let us first prove that the finite dimensional
distributions of $n^{-3/4}(X_{[nt]})_{t\ge 0}$ converge to those
of $(m\Delta_{t})_{t\ge 0}$ as $n\rightarrow\infty$. We can
rewrite for every $n\in\N$, $X_{n}=X_{n}^{(1)}+X_{n}^{(2)}$ where
$$X_{n}^{(1)}=\sum_{y\in\Z} \epsilon_{y}\Big(\sum_{i=1}^{\eta_{n-1}(y)}\xi_{i}^{(y)} -
m\Big) \; \; \; , \; \; \;  X_{n}^{(2)}=m\sum_{y\in\Z}
\eps_{y}\eta_{n-1}(y).$$ Thanks to Theorem \ref{thm11} the finite
dimensional distributions of $n^{-3/4}(X_{[nt]}^{(2)})_{t\ge 0}$
converge to those of $(m\Delta_{t})_{t\ge 0}$ as
$n\rightarrow\infty$. To conclude we show that the sequence of
random variables $n^{-3/4} (X_{n}^{(1)})_{n\in\N}$ converges for
the $L^2$-norm to 0 as $n\rightarrow +\infty$. We have
$$
\E\Big[(X_{n}^{(1)})^{2}\Big]=\E\Big[\sum_{x,y\in\Z}\eps_{x}\eps_{y}
\sum_{i=1}^{\eta_{n-1}(x)}\sum_{j=1}^{\eta_{n-1}(y)}\E[(\xi_{i}^{(x)}-m)(\xi_{j}^{(y)}-m)|\mathcal{F}
\vee \mathcal{G}]\Big]
$$
{}From the equality
$$\E[(\xi_{i}^{(x)}-m)(\xi_{j}^{(y)}-m)|\mathcal{F} \vee
\mathcal{G}]=m^2\delta_{i,j}\delta_{x,y},$$
we obtain
$$
n^{-3/2}\E\Big[(X_{n}^{(1)})^{2}\Big]=m^2 n^{-3/2}\sum_{x\in\Z}
\eta_{n-1}(x)=m^2n^{-1/2}=o(1). \; \; \; \diamond
$$

Let us recall that $M_{T_{n}}=(X_{n},Y_{n})$ for every $n\ge 1$.
The sequence of random processes $n^{-3/4}(Y_{[nt]})_{t\ge 0}$
weakly converges in ${\cal D}([0,\infty[,\R)$ to 0, thus the
sequence of $\R^2-$valued random processes
$n^{-3/4}(M_{T_{[nt]}})_{t\ge 0}$ weakly converges in ${\cal
D}([0,\infty[,\R^2)$ to the process $(m\Delta_{t},0)_{t\ge 0}$.
Theorem \ref{thm2} follows from this remark and Lemma \ref{tn}.

\section{Examples}

The main motivation of this work is the generalization of the
transience of the i.i.d. case of \cite{CP} to dependent or
inhomogeneous orientations. We obtain various extensions
corresponding to well known examples of dynamical systems such
that Bernoulli and Markov shifts,  SRB measures, rotations on the
torus, etc., our framework is very general from this point of
view. To get the transience of the walk, we need to generate the
orientations by choosing a suitable function $f$ satisfying
(\ref{C}), which requires in some sense the model not to be too
close to the deterministic case:  because to satisfy it, $f$
should not be "$\mu$-too often" 0 or 1. We describe now the
examples providing extensions of the i.i.d. case to various
disordered orientations.

{\bf 1. Shifts:} Bernoulli and Markov shifts provide the more
natural field of application of Theorem \ref{thm1},
 including  a dynamical construction of the i.i.d. case of \cite{CP} and a straightforward extension to inhomogeneous or dependent
orientations. Consider the {\em shift transformation} $T$ on the
product space $E=[0,1]^\Z$
 endowed with the Borel $\sigma$-algebra,  defined by
\begin{eqnarray*}
T: E & \longrightarrow & E\\
x=(x_{y})_{y \in \Z} & \longmapsto & (Tx)_y=x_{y+1}, \forall y\in \Z.
\end{eqnarray*}
{\em Bernoulli shifts} are considered when one starts from the
product Lebesgue measure $\mu=\lambda^{\otimes \Z}$ of the
Lebesgue measure $\lambda$ on $[0,1]$. It is $T$-invariant and
 we choose as generating function $f$  the projection on the zero coordinate:
\begin{eqnarray*}
f:E & \longrightarrow & [0,1]\\
\;x & \longmapsto & x_0.
\end{eqnarray*}
For all $y \in \Z$, we then have $f \circ T^y(x)=x_y := \xi(y) \in
[0,1]$. We consider this $\xi$'s as new random variables on $E$
whose independence is inherited from the product structure of
$\mu$. The sufficient condition (\ref{C}) becomes
\[
\int_0^1 \frac{d\lambda(x)}{\sqrt{x(1-x)}} <\infty
\] and the transience holds in this particular case. In the annealed case, the
product form of $\mu$ allows another description of the i.i.d.
case of \cite{CP}, for which we check $\xi(y)\equiv \frac{1}{2}$
for all $y\in\Z$ and $\rm{Cov}_{\mu}[\eps_0, \eps_y]=\E_{\mu}[\eps_0\eps_y]=4 \E[\xi(0) \xi(y)]-1=0$. The
result is also valid in the quenched case, for which the
distribution of the orientation has an inhomogeneous product form.

If one considers a measure $\mu$ with correlations, then the same
holds for $\pee_\mu$. Consider e.g.  $\mu$ to be a
(shift-invariant) {\em Markovian measure} on  $[0,1]^\Z$ whose
correlations are inherited from the shift via (\ref{correlation}),
with a stationary distribution $\pi$. The transience of the simple
random walk on this particular dynamically oriented lattice holds
then for $\pee_{\mu}$-a.e. environment as soon as
\[
\int_0^1 \frac{d\pi(x)}{\sqrt{x(1-x)}} < \infty.
\]

It is the case when the invariant measure $\mu$ is the usual
Lebesgue measure or Lebesgue measure of index $p$. In the quenched
case, there are no correlations by construction and the law of the
orientations depends on the measurable transformation only. This
case is nevertheless different from this of the Bernoulli shift
because the typical set of points $x$ for which the transience
holds depends on the measure $\mu$.\\

{\bf 2. SRB measures:} They provide another source of examples for
dependent orientations generated by transformations on the
interval $E=[0,1]$. A measure $\mu$ of the dynamical system $S$ is
said to be an {\em SRB} measure if the empirical measure
$\frac{1}{n} \sum_{i=1}^{n} \delta_{T^i(x)}$ converge weakly to
$\mu$ for Lebesgue a.e. $x$. There exist many other definitions of
SRB measures, see e.g. \cite{J}. In particular, it has the {\em
Bowen boundedness property} in the sense that it is close to a
Gibbs measure on some increasing cylinder, i.e. there exists a
constant $C>0$ such that for all $x \in [0,1]$ and every $n \ge 1$
\[
\frac{1}{C} \leq \frac{\mu(I_{i_1,...,i_n}(x))}{\exp{(\sum_{k=0}^{n-1}
\Phi(T^k(x)))}} \leq C
\]
where $\Phi=- \log |T'|$ and $I_{i_1,...,i_n}$ is the interval of
monotonicity for $T^n$ which contains $x$.

In some cases, it is possible to control the correlations for SRB
measures and we detail now an example where our transience result
holds, the {\em Manneville-Pomeau maps} introduced in the 1980's
to study intermittency phenomenon in the study of turbulence in
chaotic systems \cite{BPV}. They are expanding interval maps on
$E=[0,1]$ and  the original MP map is given by
\begin{eqnarray*}
T: [0,1]& \longrightarrow & [0,1]\\
x & \longmapsto & T(x)=x + x^{1+\alpha} \; {\rm mod} \; 1.
\end{eqnarray*}

The existence of an absolutely continuous (w.r.t. the Lebesgue
measure on $[0,1]$) SRB invariant measure $\mu$ has been
established by \cite{P}  and the following bounds of Radon-Nikodym
derivative $h=\frac{d\mu}{d\lambda}$ has been proved \cite{MRTMV}:
\begin{equation}\label{h}
\exists C_\star,C^\star >0 \; \rm{s.t.} \;
\frac{C_\star}{x^\alpha} < h(x) < \frac{C^\star}{x^\alpha}.
\end{equation}
This measure is known to be mixing, and a polynomial decay of
correlation, with a power $\beta >0$, has even been proved for $g$
regular enough \cite{hu,L,MRTMV,Y} :
\[
\mid C_\mu^g(y) \mid = \mathcal{O} \big(\mid y \mid ^{-\beta}
\big).
\]
The map $T$ is not invertible but we use the remark following
Theorem \ref{thm1}. It remains to find suitable function $f$ who
generates orientations for which the simple random walk is
transient. By (\ref{h}), a sufficient condition for the condition
$(\ref{C})$ to hold is
\[ \int_0^1 \frac{dx}{x^\alpha \sqrt{f(x)(1-f(x))}} < \infty
\]
and  this is for example true for the  function
$f(x)=\frac{1}{2}(1+x- T(x))$ and the choice of an $\alpha
<\frac{1}{3}$.\\

{\bf 3. Rotations:} We consider the dynamical system $S=([0,1],
{\cal B}([0,1]), \lambda, T_{\alpha})$ where $T_{\alpha}$ is the
rotation on the torus $[0,1]$ with angle $\alpha\in \R$ defined by
$$x\longmapsto x+\alpha \ \rm{ mod   }\ 1$$
and $\lambda$ is the Lebesgue measure on $[0,1]$. For every function $f:[0,1]\mapsto [0,1]$ such that $\int_0^1 f(x) \ dx=\frac{1}{2}$ and
$$ \int_0^1 \frac{dx}{\sqrt{f(x)(1-f(x))}}<\infty,$$
conclusions of Theorem \ref{thm1} hold uniformly in $\alpha$. Such
functions are called {\it admissible}. Every function uniformly
bounded from 0 and 1, with integral $\frac{1}{2}$ is admissible.
We also allow functions $f$ to take values 0 and 1: for instance,
$f_{1}(x)=x$ is admissible although $f_{2}(x)=\cos^2(2\pi x)$ is
not. We actually have no explanations about this phenomenon,
moreover we do not know the behavior (recurrence or transience) of
 the simple random walk on the dynamically oriented lattice generated by
 $f_{2}$.\\

When the generation function $f$ does not satisfy the condition
(\ref{C}), a variety of results can arise by tuning the angle
$\alpha$ to get different types of dynamical systems. Consider
$f_3={\bf 1}_{[0,1/2[ }$ and take $\alpha=\frac{1}{2q}$ for $q$ an
integer larger or equal to 1; the lattice we obtain is $\Z^2$ with
undirected vertical lines and horizontal strips of height $q$,
alternatively oriented to the left then to the right. The simple
random walk on this deterministic and periodic lattice is known to
be recurrent \cite{CP} and this provides an example of a non
ergodic system where (\ref{C}) is not fulfilled and the walk
recurrent. When the period becomes infinite, i.e. for $\alpha=0$,
the rotation is just the identity and the corresponding lattice is
$\Z^2$ with undirected vertical lines and horizontal lines all
oriented to the right (resp. all to the left) when $x\in [0,1/2[$
(resp. $x\notin [0,1/2[$). The simple random walk on this lattice
is known to be transient and this gives an example of a non
ergodic system where (\ref{C}) is not fulfilled and the walk
transient. In the ergodic case, i.e. when $\alpha$ is irrational,
we suspect the behavior of the walk to exhibit a transition
according to the type of the irrational $\alpha$. When its
approximation by rational numbers via a development in continuous
fraction is considered to be good, i.e. when its type is large,
the lattice is "quasi" periodic and the walk is believed to be
recurrent. On the other hand, the walk is believed to be transient
when this approximation is bad (when the type of the irrational is
close to 1). A deeper study of this particular choice of dynamical
system, in progress, is needed to describe more precisely the
transition between recurrence and transience in terms of the type
of the irrational.

\section{Comments}
We have extended the results of \cite{CP} to non-independent or
inhomogeneous orientations. In particular, we have proved that the
simple
   random walk is still transient for a large class of models.
   As the walk can be  recurrent for deterministic orientations,
    it would be interesting to perturb deterministic cases in order to get a
    full picture of the transience versus recurrence properties
    and a more systematic study of this problem is in progress. We
    believe that the functional limit theorem could be extended, at least
    to other ergodic dynamical systems, but this requires new
    results on random walks in ergodic random sceneries. In the i.i.d. case, Campanino {\em et al.} have also proved
    an improvement of the strong law of large numbers for the random walk in the random scenery $Z_n$: almost surely,  $\frac{Z_n}{n^\beta} \longrightarrow 0$ for all
    $\beta > \frac{3}{4}$.
     Together with our functional limit theorem and the standard results for the vertical walk, this suggests the conjecture of a
      local limit theorem, getting
    a full picture of "purely random cases", for which the
    condition on the generation $f$ holds. This work is in progress
    and we also investigate the limit theorems in more general
    cases.

%{\bf Acknowledgments:} We thank D. P\'etritis for having
%introduced us to this subject.

%\medskip

%Nadine Guillotin-Plantard\\
%Universit\'e Claude Bernard Lyon I\\
 %50, av. Tony-Garnier, Domaine de Gerland\\
   %69366 Lyon Cedex 07, France\\
    % E-mail: nadine.guillotin@univ-lyon1.fr\\

     %\medskip

      %Arnaud Le Ny\\
%Universit\'e de Paris-Sud\\
%Laboratoire de math\'ematiques\\
 %b\^atiment 425\\
 %91405 Orsay cedex, France\\
 % E-mail: arnaud.leny@math.u-psud.fr

\end{document}